\documentclass{amsproc}
\usepackage{amssymb}
\usepackage{latexsym} 
\usepackage{amsfonts} 
\usepackage{amsmath}
\usepackage{eucal} 
\usepackage{bm} 
\usepackage{bbm} 
\usepackage{graphicx} 
\usepackage[english]{varioref} 
\usepackage[nice]{nicefrac} 
\usepackage[all]{xy}
\usepackage{amsthm}
\usepackage{color}

\newcommand{\Ad}{\text {\rm Ad}}

\def\ge{\geqslant}
\def\le{\leqslant}
\def\a{\alpha}

\def\d{\delta}
\def\D{\Delta}

\def\e{\epsilon}

\def\s{\sigma}

\def\l{\lambda}

\def\i{^{-1}}
\def\ZZ{\mathbb Z}
\def\NN{\mathbb N}

\newcommand{\kk}{\Bbbk}

\def\co{\mathcal O}

\def\cw{\mathcal W}

\def\tW{\tilde W}

\def\fo{\mathfrak o}

\def\<{\langle} 
\def\>{\rangle}

\theoremstyle{plain}
\newtheorem{thm}{Theorem}[section] 
\newtheorem*{thm*}{Theorem} 
 
 \newtheorem{lem}[thm]{Lemma}

\theoremstyle{definition}

\newtheorem{example}[thm]{Example}

\theoremstyle{remark}

\newtheorem*{rmk*}{Remark}
\newtheorem*{claim*}{Claim}

\begin{document}

\author{Chuying Fang}
\address{Department of Mathematics, The Hong Kong University of Science and Technology, Clear Water Bay, Kowloon, Hong Kong}
\email{macyfang@ust.hk}
\author{Xuhua He}
\address{Department of Mathematics, The Hong Kong University of Science and Technology, Clear Water Bay, Kowloon, Hong Kong}
\email{maxhhe@ust.hk}
\thanks{X.H. is partially supported by HKRGC grants 601409.}
\title[]{Notes on partial conjugation}
\keywords{Conjugacy class, $G$-stable pieces, Weyl group, loop group}
\subjclass[2000] {20F55, 20E45, 20G25}
\dedicatory{To Fiona with adoration}

\begin{abstract}
In this notes, we will give an exposition of some results on the method of partial conjugation action. We first discuss the partial conjugation action of a parabolic subgroup of a Coxeter group. We then discuss some applications to Lusztig's $G$-stable pieces and its affine generalization. We also discuss some recent work on the $\s$-conjugacy classes of loop groups and affine Deligne-Lusztig varieties.  
\end{abstract}

\maketitle
 
\section*{Introduction} 

Let $G$ be a connected reductive algebraic group over an algebraically closed field $\kk$. Let $X$ be a ``nice'' $G \times G$-variety (e.g., a wonderful compactification of an adjoint group). In \cite{L1}, Lusztig introduced a partition of $X$ into finitely many locally closed subvarieties, which are stable under the diagonal $G$-action. These subvarieties are called $G$-stable pieces. 
Each $G$-stable pieces is a union of orbits for the diagonal action of $G$ and there is a natural bijection between these $G$-orbits and twisted conjugacy classes of a Levi factor of $G$. 
 
The partition of $X$ into $G$-stable pieces is central to Lusztig's work on parabolic character sheaves and the work of Evens and Lu \cite{EL} on Poisson geometry. More recently, Pink, Wedhorn and Ziegler \cite{PWZ} used $G$-stable pieces to study some problems in arithmetic algebraic geometry. 

The notion of $G$-stable pieces is relative new and a little hard to penetrate as the statements are fairly technical. The main purpose of this article is to explain how the method of partial conjugation action introduced by the second named author in \cite{He072} can be used to understand Lusztig's $G$-stable pieces. We will also discuss some recent progress on the study of loop groups and affine Deligne-Lusztig varieties using the method of partial conjugation action. 

\section{Partial conjugation in a Coxeter group}

\subsection{} We first recall the definition of Coxeter groups. 

Let $S$ be a finite set and $(m_{i j})_{i, j \in S}$ be a matrix with entries in $\mathbb N \cup \{\infty\}$ such that $m_{i i}=1$ and $m_{i j}=m_{j i} \ge 2$ for all $i \neq j$. Let $W$ be a group defined by the generators $s_i$ for $i \in I$ and the relations $(s_i s_j)^{m_{i j}}=1$ for $i, j \in I$ with $m_{i j}< \infty$. We say that $(W, S)$ is a {\it Coxeter group}. Sometimes we just call $W$ itself a Coxeter group.

We denote by $\ell$ the length function and $\le$ the Bruhat order. For $J \subset S$, we denote by $W_J$ the standard parabolic subgroup of $W$ generated by $J$ and by $W^J$ (resp. ${}^J W$) the set of minimal coset representatives in $W/W_J$ (resp. $W_J \backslash W$). For $J, K \subset I$, we simply write $W^J \cap {}^K W$ as ${}^K W^J$.

\subsection{} Let $J, J' \subset S$ and $\d: W_J \to W_{J'}$ such that $\d(J)=J'$. Now we consider the $\d$-twisted partial conjugation of $W_J$ on $W$ defined by $$w \cdot_\d w'=\d(w) w' w \i.$$ If $J=J'$ and $\d$ is the identity map, then we simply write $\cdot$ for $\cdot_\d$. 

For any subset $V$ of $W$, we set $$W_J \cdot_\d V=\{\d(w) w' w \i; w \in W_J, w' \in V\}.$$

The classification of $W_J$-orbits on $W$ for this twisted partial conjugation action is as follows. See \cite[Section 2]{He072}. 

\begin{thm}\label{w-dec}
For any $w \in W^J$, set $$I(J, w, \d)=\max\{K \subset J; \forall k \in K, \exists k' \in \d(K), w s_k=s_{k'} w \}.$$ Then 

(1) $W=\sqcup_{w \in W^J} W_J \cdot_\d (w W_{I(J, w, \d)})$.

(2)  The inclusion $W_{I(J, w, \d)} \to W_J \cdot_\d (w W_{I(J, w, \d)})$ induces a natural bijection between the $\Ad(w) \i \circ \d$-twisted conjugacy classes on $W_{I(J, w, \d)}$ and the $W_J$-orbits on $W_J \cdot_\d (w W_{I(J, w, \d)})$. 
\end{thm}

We'll say a few words on the $\Ad(w) \i \circ \d$-twisted conjugation action in (2) of the above theorem. It is easy to see from the definition of $I(J, w, \d)$ that for any $v \in W_{I(J, w, \d)}$, $w \i \d(v) w \in W_{I(J, w, \d)}$. The $\Ad(w) \i \circ \d$-twisted conjugation action of $W_{I(J, w, \d)}$ on itself is defined to be $v \cdot_{\Ad(w) \i \circ \d} v'=w \i \d(v) w v' v \i \in W_{I(J, w, \d)}$. Moreover, for any $v, v'\in W_{I(J, w, \d)}$, $$w (v \cdot_{\Ad(w) \i \circ \d} v')=w (w \i \d(v) w v' v \i)=\d(v) w v' v \i =v \cdot_\d (w v').$$ 

\subsection{} In particular, we have a map $$\pi_{J, \d}: W \to W^J, W_J \cdot_\d (w W_{I(J, w, \d)}) \mapsto w$$ that is constant on each orbit for the $\d$-twisted partial conjugation action of $W_J$. 

Below are two examples. Here red color for the elements in $W^J$ and blue color for the elements in $W_J$. 

\begin{example} $W=S_4$, $J=J'=\{1, 2\}$ and $\d$ is the identity map. 

The image of $w=s_2 s_3 s_2 s_1 s_2$ is obtained as follows. 
\begin{align*} \textcolor{red}{s_2 s_3} \textcolor{blue}{s_2 s_1 s_2} & \xrightarrow{\text{ conj. by } s_2}  \textcolor{red}{s_3}  \textcolor{blue}{s_2 s_1} \xrightarrow{\text{ conj. by } s_1} \textcolor{red}{s_3} \textcolor{blue}{s_1 s_2} \xrightarrow{\text{ conj. by } s_2} \textcolor{red}{s_2 s_3} \textcolor{blue}{s_1}  \\ & \xrightarrow{\text{ conj. by } s_1}\textcolor{red}{s_1 s_2 s_3} \in W^J.
\end{align*}

So $\pi_{J, \d} (s_2 s_3 s_2 s_1 s_2)=s_1 s_2 s_3$. 

\end{example}

\begin{example} $W=S_4$, $J=J'=\{1, 3\}$ and $\d$ is the identity map.

The image of $w=s_2 s_1 s_3 s_2 s_1$ is obtained as follows. 
\begin{align*} \textcolor{red}{s_2 s_1 s_3 s_2}  \textcolor{blue}{s_1} \xrightarrow{\text{ conj. by } s_1} \textcolor{red}{s_2 s_1 s_3 s_2} \textcolor{blue}{s_3} \xrightarrow{\text{ conj. by } s_3} \textcolor{red}{s_2 s_1 s_3 s_2}  \textcolor{blue}{s_1}.
\end{align*}

In this case, $I(J, s_2 s_1 s_3 s_2, \d)=\{1, 3\}$ and $\pi_{J, \d} (s_2 s_1 s_3 s_2 s_1)=s_2 s_1 s_3 s_2$. 
\end{example}

\subsection{} In order to use the partial conjugation of Weyl group to understand structure of reductive groups, we need additional properties related to the partial conjugation. 

The following notion is motivated by \cite{GP93}.

Let $w, w' \in W$. We write $w \xrightarrow{i}_\d w'$ for $i \in J$ if $w'=s_{\d(i)} w s_i$ and $\ell(w') \le \ell(w)$. We write $w \to_{J, \d} w'$ if there is a sequence $w=w_1 \xrightarrow{i_1}_\d w_2 \xrightarrow{i_2}_\d \cdots \xrightarrow{i_n}_\d w_{n+1}=w'$, where $i_1, \cdots, i_n \in J$. We write $w \approx_{J, \d} w'$ if $w \to_{J, \d} w'$ and $w' \to_{J, \d} w$. It is easy to see that $w \approx_{J, \d} w'$ if and only if $w \to_{J, \d} w'$ and $\ell(w)=\ell(w')$. If $J=S$ and $\d$ is the identity map, then we simply write $\to$ for $\to_{J, \d}$ and $\approx$ for $\approx_{J, \d}$. 

We say that $w \in W$ is terminal with respect to $(J, \d)$ if for any $w' \in W$ with $w \to_{J, \d} w'$, we have that $w \approx_{J, \d} w'$. 

Now we have the following results \cite[Section 3]{He072}.

\begin{thm}\label{min}
Let $w \in W$ with $\pi_{J, \d}=w'$. Then there exists $x \in W_{I(J, w, \d)}$ such that $w \to_{J, \d} w' x$. If moreover, $\ell(w)=\ell(w')$, then $w \approx_{J, \d} w'$. 
\end{thm}

In fact, if $W_J$ is a finite Coxeter group, then we can choose $x$ in such a way that $w' x$ is a minimal length element in the $\d$-twisted conjugacy class of $W_J$ that contains $w$. 

\subsection{} By \cite[Corollary 4.5]{He072}, for any $\d$-twisted $W_J$-conjugacy class $\co$ that contains an element in $W^J$, the following conditions are equivalent:

(1) $v$ is a minimal element in $\co$ with respect to the restriction to $\co$ of the Bruhat order on $W$.

(2) $v$ is an element of minimal length in $\co$.

We denote by $\co_{\min}$ the set of elements in $\co$ that satisfy the above conditions. As in \cite[4.7]{He072}, we have a natural partial order $\le_{J, \d}$ on $W^J$ defined as follows:

Let $w, w' \in W^J$. Then $w \le_{J, \d} w'$ if for some (or equivalently, any) $v' \in (W_J \cdot_\d w')_{\min}$, there exists $v \in (W_J \cdot_\d w)_{\min}$ such that $v \le v'$.

In general, for $w \in W^J$ and $w' \in W$, we write $w \le_{J, \d} w'$ if there exists $v \in (W_J \cdot_\d w)_{\min}$ such that $v \le w'$.

\section{Partial conjugation in groups with $(B, N)$-pair}

\subsection{} We first recall the definition of a $(B, N)$-pair. 

A group $G$ has a $(B, N)$-pair if the following axioms hold:

\begin{enumerate}

\item $G$ is generated by the subgroups $B$ and $N$.

\item $B \cap N$ is a normal subgroup of $N$.

\item The group $W=N/B \cap N$ is generated by a set of elements $s_i$ of order $2$, for $i$ in some nonempty set $S$.

\item For any $i \in S$ and $w \in W$, $s_i B w \subset B s_i w B \cup B w B$. 

\item $s_i \notin N_G(B)$ for all $i \in S$. 
\end{enumerate}

In this case, $(W, S)$ is a Coxeter group and $G=\sqcup_{w \in W} B w B$. Moreover, for any $i \in S$ and $w \in W$, we have that \[B s_i B w B=\begin{cases} B s_i w B, & \text{ if } s_i w>w; \\ B s_i w B \sqcup B w B, & \text{ if } s_i w<w.\end{cases}\] \[B w B s_i B=\begin{cases} B w s_i B, & \text{ if } w s_i>w; \\ B w s_i B \sqcup B w B, & \text{ if } w s_i<w.\end{cases}\]

\subsection{} Let $J, J' \subset S$ and $L_J, L_{J'}$ the corresponding Levi subgroups of $G$. Let $\s: L_J \to L_{J'}$ a isomorphism of (abstract) groups that sends $B \cap L_J$ to $B \cap L_{J'}$ and $H \cap L_J$ to $H \cap L_{J'}$. Then $\s$ induces an isomorphism $\d: W_J \to W_{J'}$ such that $\d(J)=J'$. Define the $\s$-twisted partial conjugation action of $L_J$ on $G$ by $l \cdot_\s g=\s(l) g l \i$. 

The following Lemma is crucial in the study of $\s$-twisted partial conjugation action. 

\begin{lem}\label{1}
We keep the notations as above. Let $w, w' \in W$. 

(1) If $w \approx_{J, \d} w'$, then $L_J \cdot_\s B w B=L_J \cdot_\s B w' B$. 

(2) If $w$ is not terminal with respect to $(J, \d)$, then $$L_J \cdot_\s B w B \subset \cup_{v \in W, \ell(v)<\ell(w)} L_J \cdot_\s B v B.$$
\end{lem}

\

Now similar to the decomposition in Theorem \ref{w-dec}, we have the following result. The proof is based on combinatorial properties of partial conjugation in Weyl group. 

\begin{thm}\label{dec-0}
$G=\cup_{w \in W^J, x \in W_{I(J, w, \d)}} L_J \cdot_\s B w x B$. 
\end{thm}

Proof. Notice that $G=\sqcup_{v \in W} B v B$. Thus it suffices to prove that for any $v$, \[\tag{1} B v B \subset \cup_{w \in W^J, x \in W_{I(J, w, \d)}} L_J \cdot_\s B w x B.\]

We argue by induction on $\ell(v)$. Let $w_1=\pi_{J, \d}(v)$. Then by Theorem \ref{min}, there exists $x_1 \in W_{I(J, w_1, \d)}$ such that $v \to_{J, \d} w_1 x_1$. Now by Lemma \ref{1}, \[\tag{2} B v B \subset L_J \cdot_\s B w_1 x_1 B \cup \cup_{v' \in W, \ell(v')<\ell(v)} L_J \cdot_\s B v' B.\] By induction hypothesis, for any $v'$ with $\ell(v')<\ell(v)$, \[\tag{3} B v' B \subset \cup_{w \in W^J, x \in W_{I(J, w, \d)}} L_J \cdot_\s B w x B.\] Now (1) follows from (2) and (3). This finishes the proof. \qed

\

It is very interesting to consider when $J$ is of finite type. We discuss in the following two sections the case where $G$ is of finite type and affine type. These results can also be generalized to arbitrary Kac-Moody groups. We do not go into this in the article. 

\section{Partial conjugation: finite type}

The notion of $G$-stable pieces was introduced by Lusztig in \cite{L1} in the study of parabolic character sheaves. A simpler formulation was given in \cite{He06} and the closure relation was found in \cite{He07}. Different approaches were obtained by Evens and Lu in \cite{EL} and by Springer in \cite{Sp3}. The notion was later generalized by Lu and Yakimov to $R$-stable pieces in \cite{LY}.  We refer to the survey article of Springer \cite{Sp1} for some applications of $G$-stable pieces in representation theory. 

\subsection{}\label{3.1} In this section, assume that $G$ is a reductive group over an algebraically closed field $\kk$. Let $B$ be a Borel subgroup of $G$ and $T \subset B$ a maximal torus. We denote by $W$ the Weyl group of $G$ and $S$ the set of simple roots. For any subset $J$ of $W$, we denote by $w_0^J$ the maximal element in $W_J$. 

Assume that $J \subset S$ and $\s: L_J \to L_{J'}$ is an isomorphism of algebraic groups. Then for any $w \in W^J$, the map $L_{I(J, w, \d)} \to L_{I(J, w, \d)}$, $l \mapsto w \i \s(l) w$ is again an isomorphism of algebraic groups. Therefore the map \[L_{I(J, w, \d)} \times (B \cap L_{I(J, w, \d)}) \to L_{I(J, w, \d)}, \quad (l, b) \mapsto w \i \s(l) w b l \i\] is surjective (see \cite{Ste}). As a consequence, one can see that $$L_{I(J, w, \d)} \cdot_\s B w x B \subset L_{I(J, w, \d)} B w B=L_{I(J, w, \d)} \cdot_\s B w B$$ for all $x \in W_{I(J, w, \d)}$. Now by Theorem \ref{dec-0}, $G=\cup_{w \in W^J} L_J \cdot_\s B w B$. In fact, a detailed analysis shows that this union is a disjoint union (see \cite{L1}, \cite{Sp3}). 

\begin{thm}
$G=\sqcup_{w \in W^J} L_J \cdot_\s B w B$.
\end{thm}

Now we describe the closure relations between these subvarieties of $G$. The proof is based on combinatorial property of partial conjugation in Weyl group. Details can be found in \cite[Section 5]{He072}.

\begin{thm}
For any $w \in W$, $$\overline{L_J \cdot_\s B w B}=\sqcup_{w' \in W^J, w' \le_{J, \d} w} L_J \cdot_\s B w' B.$$ 
\end{thm}







\subsection{} We set $R_{J, \s}=\{(l u, \s(l) u'); l \in L_J, u \in U_{P_J}, u' \in U_{P_{\d(J)}}\} \subset G \times G$, where $U_{P_J}$ is the unipotent radical of the standard parabolic subgroup $P_J$ and $U_{P_{\d(J)}}$ is the unipotent radical of $P_{\d(J)}$. The action of $R_{J, \d}$ on $G$ is defined by $(l u, \s(l) u', g) \mapsto l u g (\s(l) u') \i$. Then it is easy to see that for any $w \in W$, $R_{J, \s} \cdot B w B=(L_J \cdot_\d B w B) \i$. Thus 

(1) $G=\sqcup_{w \in W^J} R_{J, \s} \cdot B w \i B$.

(2) For any $w \in W^J$, $\overline{R_{J, \s} \cdot B w B}=\sqcup_{w' \in W^J, w' \le_{J, \d} w} R_{J, \s} \cdot B (w') \i B$. 

\subsection{} Let $G_\D$ be the diagonal image of $G$ in $G \times G$. Then we may regard $G$ as $G \times G/G_\D$ via $g \mapsto (1, g) G_\D/G_\D$ and thus there exists a natural bijection of the $R_{J, \s}$-orbits on $G \cong (G \times G)/G_\D$ and $G_\D$-orbits on $(G \times G)/R_{J, \s}$. This bijection is good in the sense of \cite[Lemma 1.6]{Sp3}. Under this bijection, the subset $R_{J, \d} \cdot B w \i B$ of $G$ corresponds to the subset $G_\D \cdot (B w, B) R_{J, \s}/R_{J, \s}$ of $(G \times G)/R_{J, \s}$. For any $w \in W^J$, we write $Z_{J, w, \s}$ for $G_\D \cdot (B w, B) R_{J, \s}/R_{J, \s}$ and call it a {\it $G$-stable piece} of $(G \times G)/R_{J, \s}$. 

Below are some properties of the $G$-stable pieces. 

(1) $(G \times G)/R_{J, \s}=\sqcup_{w \in W^J} Z_{J, w, \s}$. See \cite{L1}.

(2) For any $w \in W^J$, $Z_{J, w, \s}$ is a locally closed smooth subvariety of $(G \times G)/R_{J, \s}$ of codimension $\ell(w_0^J w_0)-\ell(w)$. See \cite{L1}. 

(3) There is a natural bijection between the $G_\D$-orbits on $Z_{J, w, \s}$ and the $\s \circ \Ad(w)$-twisted conjugacy classes of $L_{I(J, w, \d)}$. See \cite{L1}.

(4) For any $w \in W^J$, $\overline{Z_{J, w, \s}}=\sqcup_{w' \in W^J, w' \le_{J, \d} w} Z_{J, w', \s}$. See \cite{He072}.

(5) For any $w \in W^J$, $\overline{Z_{J, w, \s}}$ is a semi-normal variety. See \cite{HT}.

\section{Some variations}

\subsection{} We keep the notations in Section 3.

The action of $-w_0$ on the set of simple roots gives an involution $*: S \to S$. Let $J \subset S$ and $P^-_J$ be the parabolic subgroup of $G$ opposite to the standard parabolic $P_J$. Then $P^-_J={}^{w_0^J w_0} P_{J^*}$. Define $$Z_J=(G \times G) \times_{P^-_J \times P_J} L_J,$$ here the action of $P^-_J  \times P_J$ on $G \times G \times L_J$ is defined by $(p, p') \cdot (g, g', l)=(g p \i, g' (p') \i, \pi_{P_J^-}(p) l \pi_{P_J}(p'))$, where $\pi_{P_J^-}: P^-_J \to L_J$ and $\pi_{P_J}: P_J \to L_J$ are projection maps. It is easy to see that $Z_J \cong (G \times G)/R^-_J$, where $R^-_J=\{(l u, l u'); l \in L_J, u \in U_{P^-_J}, u' \in U_{P_J}\}$. 

Now define $\s: L_{J^*} \to L_J$ by $\s(l)=(w_0^J w_0) l (w_0^J w_0) \i$. Then $$R^-=(w_0^J w_0, 1) R_{J^*, \s} (w_0^J w_0, 1) \i.$$ Hence the isomorphism $G \times G \to G \times G$, $(g, g') \mapsto (g (w_0^J w_0) \i, g')$ induces an isomorphism $G \times G/R_{J^*, \s} \to Z_J$. 

Now for any $w \in W^J$, we write $Z_{J, w}$ for $G_\D \cdot (B w, B) R^-_J/R^-_J$ and call it a {\it $G$-stable piece} of $Z_J$.  

\subsection{} In this subsection, we assume that $G$ is of adjoint type. In the same manner as $Z_J$, we define $X_J=(G \times G)_{P_J^- \times P_J} L_J/Z(L_J)$. This is a boundary $G \times G$-orbit of the wonderful compactification $X$ of $G$, see \cite{DP}.

We denote by $h_J$ the image of $(1, 1, 1)$ in $X_J$. For $w \in W^J$, set $$X_{J, w}=G_\D \cdot (B w, B) h_J$$ and call it a $G$-stable piece of $X$. Below are some properties of $X_{J, w}$. 

(1) $X=\sqcup_{J \subset S, w \in W^J} X_{J, w}$. See \cite{L1}.

(2) For any $J \subset S, w \in W^J$, $X_{J, w}$ is a locally closed smooth subvariety of $X$ of codimension $\ell(w)+\sharp(S-J)$. See \cite{L1}.

(3) For any $J \subset S, w \in W^J$, $\overline{X_{J, w}}=\sqcup_{K \subset J, w' \in W^K, w \le_{J, id} w'} X_{K, w'}$. See \cite{He07}. 

(4) For any $J \subset S, w \in W^J$, $\overline{X_{J, w}}$ is a semi-normal variety. See \cite{HT}.

(5) If $\overline{X_{J, w}}$ contains finitely many $G_\D$-orbits, then it admits a cellular decomposition. See \cite{He07}. 

\subsection{} In the rest of this section, we assume that the characteristic of the ground field $\kk$ is positive. Instead of considering isomorphism of algebraic groups $\s: L_J \to L_{J'}$ in $\S$\ref{3.1}, we may consider the morphism $\s \circ F: L_J \to L_{J'}$, where $F: L_J \to L_J$ is a Frobenius morphism. Now using Lang's theorem, one can show that for any $w \in W^J$, $R_{J, \s \circ F}$ acts transitively on $R_{J, \s \circ F} \cdot B w \i B$. 

Similarly, for any Frobenius morphism $F: G \to G$ and any $w \in W^J$, $G_F$ acts transitively on $G_F \cdot (B w, B) R_{J, \s} \subset G \times G/R_{J, \s}$, where $G_F=\{(g, F(g)); g \in G\} \subset G \times G$. The closure relation between the $G_F$-orbits is still given by $\le_{J, \d}$. This is used in the study of finer Deligne-Lusztig varieties in a partial flag variety \cite{He09}. This is also used by Pink, Wedhorn and Ziegler, in the study of $F$-zips \cite{PWZ}. 

\subsection{} Now we assume that $J=J_1 \sqcup J_2$ and $W_J=W_{J_1} \times W_{J_2}$. We consider the actions of $W_{J_2} \times W_{\d(J_2)}$ on the set of $W_{J_1}$-conjugacy classes of $W$ and its applications to reductive groups. We only list the results below. The details will appear in a future work. 

For $w \in W^{J_1}$, define $K_i(w)$ (for $i \ge 0$) as follows: $K_0(w)=J_2$ and $K_i(w)=w K_{i-1}(w) \cap J_1$ for $i \ge 1$. Set $$\cw(J_1, J_2)=\{w \in W^{J_1}; w \in {}^{K_i(w)} W \text{ for all } i\}.$$

This set first appeared in the study of conjugacy classes in reductive monoid. See \cite{P}.

We have that

(1) $W=\sqcup_{w \in \cw(J_1, J_2)} W_{J_2} \pi_{J_1, id} \i(w) W_{J_2}$.

(2) For any $w \in \cw(J_1, J_2, \d)$, the $(W_{\d(J_2)} \times W_{J_2}) (W_{J_1})_{id}$-orbits on $W_{J_2} \pi_{J_1, id} \i(w) W_{J_2}$ are in bijection with the $(W_{J_1})_{id}$-orbits on $\pi_{J_1, \d} \i (w)$, i.e. $\Ad(w)$-twisted conjugacy classes on $W_{I(J_1, w, id)}$. 

Now let $R_{J_1, J_2}=\{(l p_1, l p_2; l \in L_{J_1}, p_1, p_2 \in U_{P_{J_1 \cup J_2}} L_{J_2}\}$. Then $$(G \times G)/R_{J_1, J_2}=\sqcup_{w \in \cw(J_1, J_2)} G_F \cdot (B w, B) R_{J_1, J_2}/R_{J_1, J_2}$$ and for any $w \in \cw(J_1, J_2)$, $G_F \cdot (B w, B) R_{J_1, J_2}/R_{J_1, J_2}$ is a single $G_F$-orbit. 

Notice that if $G=GL_n$, then each $G \times G$-orbit on $\mathfrak{gl}_n$ is of the form $G \times G/R_{J_1, J_2}$, where $J_1=\{1, 2, \cdots, i-1\}$ and $J_2=\{i+1, i+2, \cdots, n-1\}$ for some $i$. In particular, there are finitely many $G_F$-orbits on $\mathfrak g \mathfrak l_n$. These orbits classify the isomorphism classes of restricted abelian Lie algebras of dimension $n$. Thus we obtain the following unpublished result of Lin \cite{Lin}.

\begin{thm}
There are only finitely many isomorphism classes of restricted Lie algebra structure on an abelian Lie algebra of dimension $n$. 
\end{thm}

\section{Partial conjugation: affine type} 

\subsection{} Let $L=\kk((\e))$ be the formal Laurent series and $\fo=\kk[[\e]]$ be the ring of formal power series. Let $K=G(\fo)$ be a maximal bounded subgroup of the loop group $G(L)$. Let $I$ be the inverse image of $B^-$ under the projection map $K \mapsto G(\kk)$ sending $\e$ to $0$ and $K_1$ be the kernel of this projection map. Let $\tW=N(T(L))/(T(L) \cap I)$ be the extended affine Weyl group of $G(L)$. It is known that $\tW=W \ltimes X_*(T)=\{w \e^\chi; w \in W, \chi \in X_*(T)\}$, where $X_*(T)$ is the coweight lattice. 

We only discuss in the section the partial conjugation of $K$ on $G(L)$. The partial conjugation of arbitrary parahoric subgroup on $G(L)$ can be found in \cite{L2}. 

Let $\tilde W^S$ be the set of minimal coset representatives in $\tilde W/W$. Define $K_w=K \cdot I w I$ for $w \in \tilde W^S$ and call it a {\it $K$-stable piece} of $G(L)$. Here $\cdot$ means conjugation action. Similar to our discussion in $\S 3$, we have that 

(1) $G(L)=\sqcup_{w \in \tilde W^S} K_w$. See \cite{L2}. 

(2) Let $R=\{(g u, g u'); g \in G(\kk), u, u' \in K_1\}$. Then there is a natural bijection between the $R$-orbits on $K_w$ and the $\Ad(w)$-twisted conjugacy class of $L_{I(S, w, id)}(\kk)$. See \cite{L2}.

(3) For any $w \in \tilde W^S$, $\overline{K_w}=\sqcup_{w' \in \tilde W^S, w' \le_{S, id} w} K_{w'}$. See \cite{He4}. 

The Frobenius-twisted case is studied by Viehmann in \cite{V2}. 

\subsection{} In this subsection, we assume that $G$ is an adjoint group. The specialization $\e \mapsto 0$ defines a map $s$ from the loop group $G(L)$ to the wonderful compactification $X$ of $G(\kk)$. This specialization map was introduced by Springer in \cite[2.1]{Sp2}. It is easy to see that $s(K)=G$ and $s(I)=B^-$.

Notice that any element in $\tilde W^S$ is of the form $x \e^{-\l}$ for some dominant coweight $\l$ and $x \in W^{I(\l)}$, where $I(\l)=\{i \in S; \<\l, \a_i\>=0\}$. Now we have the following correspondence between the $K$-stable piece in $G(L)$ and the $G$-stable piece in the wonderful compactification $X$ of $G(\kk)$. See \cite{He4}.

(1) For any $\l \in Y^+$ and $x \in W^{I(\l)}$, $s(K_{x \e^{-\l}})=X_{I(-w_0 \l), w_0 x w_0}$. 

(2) For any $J \subset S$ and $x \in W^J$, $s \i (X_{J, x})=\sqcup_{\l \in Y^+, I(\l)=-w_0 J} K_{w_0 x w_0 \e^{-\l}}$.

This correspondence is a key ingredient in \cite{He4} in the study of the relation between the closure of unipotent variety of $G(\kk)$ in $X$ and the affine Deligne-Lusztig varieties in the affine flag of the loop group $G(L)$. 

\section{Conjugation action in loop groups}

\subsection{} It is a challenging problem to study the case where $J$ is not of finite type in the setting of Section 2.  

In this section, we only consider the special case where $J=S$ is of affine type. More precisely, we consider a ``twisted'' conjugation action of $G(L)$ on itself as $g \cdot_\s h=\s(g) h g \i$ for $g, h \in G(L)$. Here $\s$ is a bijective group homomorphism on $G(L)$ of one of the following type:

(1) For any nonzero element $a \in \kk$, define $\s_a(p(\e))=p(a \cdot \e)$ for any formal Laurent power series $p(\e)$. We extend $\s_a$ to a group homomorphism on $G(L)$, which we still denote by $\s_a$.

(2) If $\kk$ is of positive characteristic and $F: \kk \to \kk$ is a Frobenius morphism. Then set $F(\sum a_n \e^n)=\sum F(a_n) \e^n$. We extend $F$ to a group homomorphism on $G(L)$, which we denote by $\s_F$. 

The $\s_a$-conjugacy classes were studies by Baranovsky and Ginzburg in \cite{BG96}. The $\s_F$-conjugacy classes were studied by Kottwitz in \cite{Ko97}. 

For simplicity, we focus on the case where $G=GL_n$. In this case, $\tilde W=S_n \ltimes \ZZ^n$. The discussion for other groups can be found in \cite{He5} and \cite{GH}.

\subsection{} Similar to our discussion in the previous sections, we will use some special properties of the affine Weyl group to understand the twisted conjugation action in loop groups. 

We call an element $w \in \tilde W$ a {\it good element} if $\ell(w^n)=n \ell(w)$ for all $n \in \NN$ and we write $\tW_{good}$ for the set of all good elements in $\tW$. 

For any $J \subset \tilde S$, set $\cw_J=\{w x; w \in \tW^J, x \in W_{I(J, w, id)}\}$. By Theorem \ref{min}, any $W_J$-conjugacy class in $\tW$ contains an element in $\cw_J$ that is a minimal length element in that conjugacy class. We call a conjugacy class $\co$ of $\tW$ a {\it distinguished conjugacy class} if $\co_{\min} \cap \cw_J \subset \tW^J$ for all $J \subsetneqq \tilde S$. 

The following result was proved in \cite{He5}. 

\begin{thm}\label{aff-min}
Let $\co$ be a conjugacy class in $\tW$. Then

(1) For any $w \in \co$, there exists $w' \in \co_{\min}$ such that $w \to w'$. 

(2) $\co$ is distinguished if and only if $\co$ contains a good element. In this case, $\co_{\min} \subset \tW_{good}$. 

(3) If $\co$ is distinguished, then for any $w, w' \in \co_{\min}$, $w \approx w'$. 
\end{thm}

\subsection{} For any distinguished conjugacy class $\co$ of $\tW$, we choose a minimal length element $w_\co$ in $\co$. By \cite[8.3]{He5}, we introduce a partial order on the set of distinguished conjugacy classes as follows. We write $\co \le \co'$ if there exists $w \in \co_{\min}$ such that $w \le w_{\co'}$. It is showed that this definition is independent of the choice of the representative $w_{\co'}$. 

In general, for any $x \in \tW$, we write $\co \le x$ if there exists $w \in \co_{\min}$ such that $w \le x$. 

\

Based on these special properties of $\tW$, we have the following result. See \cite[Section 11]{He5}. 

\begin{thm} For $\s=\s_a$ or $\s_F$, we have that

(1) $G(L)=\sqcup_{\co \text{ distinguished }} G(L) \cdot_\s I w_\co I$.

(2) For any $x \in \tW$, $\overline{G(L) \cdot_\s I x I}=\sqcup_{\co \text{ distinguished }, \co \le x} G(L) \cdot_\s I w_\co I$. 
\end{thm}

\begin{rmk*}
If $\s=\s_F$, then for any good element $w$, $G(L) \cdot_\s I w I$ is a single $\s_F$-conjugacy class of $G(L)$. In this case, part (1) of the theorem is a reformation of Kottwitz's classification of $\s_F$-conjugacy classes of $G(L)$. The closure relation of is $\s_F$-conjugacy classes was also obtained by Viehmann in \cite{V2}. 
\end{rmk*}

\subsection{} If $\s=\s_F$, then each $\s$-conjugacy class of $G(L)$ contains a representative $b$, here $b$ is a block diagonal matrix and each block is of the form $$\left( \begin{array}{cc} 0 & \e^{k_i+1} I_{k_i'} \\ \e^{k_i} I_{n_i-k'_i}  & 0 \end{array} \right).$$ The image $\underline b$ of $b$ in $\tilde W$ is in general, not a good element. However,  the $S_n$-conjugacy class of $\underline b$ contains a unique element in $\tW^S$. That element is a good element associated to the $\s$-conjugacy class of $b$. 

\subsection{} In the rest of this article, we discuss some recent progress on the affine Deligne-Lusztig varieties. 

By definition, the affine Deligne-Lusztig varieties associated to an element $w$ in the extended affine Weyl group $\tilde W$ and a base point $b \in G(L)$ is $$X_w(b)=\{g I \in G(L)/I; g \i b \s_F(g) \in I x I\}.$$ They play an important role in the study of Shimura varieties. More details can be found in the survey papers of G\"ortz \cite{Go}, Haines \cite{Ha} and Rapoport \cite{Ra}.

\subsection{} 
Define

\begin{itemize}
\item $\eta_1: \tilde W=W \ltimes X_*(T)  \to W$, projection map.

\item $\eta_2: \tilde W \to W$ such that $\eta_2(x) \i x$ lies in the dominant Weyl chamber.

\item $\eta(x)=\eta_2(x) \i \eta_1(x) \eta_2(x)$.
\end{itemize}

Then we have the following result describing the dimension of $X_w(b)$. 

\begin{thm}
Let $b \in G(L)$ be a basic element, i.e., a length $0$ element in $\tW$. Let $w$ be an element in the lowest two-sided cell of $\tW$. Then $X_w(b) \neq \emptyset$ if and only if $b$ and $w$ are in the same component of $G(L)$ and $\eta(w)$ is not in any proper standard parabolic subgroup of $W$. If moreover, the translation part of $w$ is regular, then $$\dim X_w(b)=\frac{1}{2}(\ell(w)+\ell(\eta(w))-\text{def}(b)), $$ here $\text{def}(b)=rank(G)-rank(J_b)$, where $J_b$ is the $\s_F$-centralizer of $b$. 
\end{thm}

The ``only if'' part of the theorem was proved by G\"ortz, Haines, Kottwitz and Reuman in \cite{GHKR2} and the remaining part was proved in \cite{GH}. In loc. cit., the upper bound $\dim X_w(b) \le \frac{1}{2}(\ell(w)+\ell(\eta(w))-\text{def}(b)$ was deduced from the dimension formula in affine Grassmannian \cite{GHKR1} and \cite{V}. Theorem \ref{aff-min} is a key ingredient to establish the lower bound  $\dim X_w(b) \ge \frac{1}{2}(\ell(w)+\ell(\eta(w))-\text{def}(b)$. 

\subsection{} Notice that $X_w(b)=\lim\limits_{\rightarrow} X_i$ for some closed subschemes $X_1 \subset X_2 \subset \cdots \subset X_w(b)$ of finite type. Thus we may define the Borel-Moore homology of $X_w(b)$ as $$H^{BM}_j(X_w(b), \bar {\mathbb Q}_l)=\lim\limits_{\rightarrow} H^j_c(X_i, \bar {\mathbb Q}_l)^*,$$ here  $l$ is a prime number not equal to the characteristic of $\kk$. The action of the $\s_F$-centralizer $J_b$ of $b$ on $X_w(b)$ induces a smooth representation of $J_b$ on $H^{BM}_j(X_w(b), \bar {\mathbb Q}_l)$. The following property was obtained in \cite{He5}. The proof is again based on Theorem \ref{aff-min}. 

\begin{thm}
Let $\mu$ a dominant regular coweight of $G=GL_n$. Then $J_{\e^\mu}=T(\mathbb F_q((\e)))$ and $T(\mathbb F_q[[\e]])$ acts trivially on $H^{BM}_j(X_w(\e^\mu), \overline{\mathbb Q_l})$ for all $w \in \tilde W$ and $j \in \mathbb Z$.  
\end{thm}

\section*{Acknowledgment} We are grateful to S. Evens for his comments on this paper.

\end{document}